\documentclass[a4paper,12pt]{amsart}

\newtheorem{dfn}{Definition}[section]
\newtheorem{thm}[dfn]{Theorem}

\newtheorem{rem}[dfn]{Remark}

\def\R{{\mathbb{R}}}
\def\N{{\mathbb{N}}}
\def\Z{{\mathbb{Z}}}

\usepackage{verbatim}
\usepackage{amssymb}
\usepackage{amsbsy}
\usepackage{amscd}
\usepackage{amsmath}
\usepackage{amsthm}
\usepackage[mathscr]{eucal}
\usepackage[all]{xy}
\usepackage{graphicx}
\usepackage{adjustbox}

\title[Discrete integrable systems and Pitman's transformation]{Discrete integrable systems\\and Pitman's transformation}

\author[D.~A.~Croydon]{David A.\ Croydon}
\address{Research Institute for Mathematical Sciences, Kyoto University, Kyoto 606-8502, Japan}
\email{croydon@kurims.kyoto-u.ac.jp}

\author[M.~Sasada]{Makiko Sasada}
\address{Graduate School of Mathematical Sciences, University of Tokyo, 3-8-1, Komaba, Meguro-ku, Tokyo, 153-8914, Japan}
\email{sasada@ms.u-tokyo.ac.jp}


\begin{document}

\begin{abstract}
We survey recent work that relates Pitman's transformation to a variety of classical integrable systems, including the box-ball system, the ultra-discrete and discrete KdV equations, and the ultra-discrete and discrete Toda lattice equations. It is explained how this connection enables the dynamics of the integrable systems to be initiated from infinite configurations, which is important in the study of invariant measures. In the special case of spatially independent and identically distributed configurations, progress on the latter topic is also reported.
\end{abstract}

\maketitle

\section{Introduction}

Let us recall the well-known Pitman's transformation, which for discrete-time paths $S: \Z_+ \to \R$ and continuous-time paths $S:\R_+ \to \R$ is defined by the mapping
\begin{equation}\label{pitman}
S \mapsto  2M_{\cdot} -S_{\cdot}
\end{equation}
where $M_x=\sup_{0 \le y \le x} S_y$. (In the continuous-time case, one needs to make an assumption on $S$ to ensure that $M$ is finite.) We will denote this transform by $T^{\Z}S:=2M-S$ and $T^{\R}S:=2M-S$ respectively.

Since its introduction in \cite{Pitman}, where it was shown to relate the one-dimensional Brownian motion and the three-dimensional Bessel process, Pitman's transformation and variations thereof have been extensively studied. Following Pitman's original work, much of this interest has been in the context of stochastic processes (see \cite{Bertoin, HMOC, Jeulin, MY0, MY, MY1, MY2, Rog, RogPit, Saisho}, for example). Moreover, there has been lively activity involving Pitman-type transformations in the areas of queuing theory and stochastic integrable systems, where it has links to models such as random polymers, random matrices, the quantum Toda lattice and the KPZ equation (see \cite{BBOC1,OC0,OC,OCY}, for example), and such transformations have also been used to define a continuous model of crystals \cite{BBOC2}. We highlight that one of the most important and well-studied adaptations of Pitman's transformation arising within this research is its exponential version, which is defined for a continuous-time path $S : \R_+ \to \R$ by the mapping
\[S \mapsto 2 M^{\int}_\cdot-S_\cdot,\]
where
\[M^{\int}_\cdot:=\frac{1}{2}\log \left(\int_0^{\cdot} \exp(2S_y) dy\right).\]
(Again, it is necessary to make suitable assumptions on $S$ to ensure that the above transformation is well-defined.) We will denote this transform by $T^{\int}S:=2M^{\int}-S$.

The purpose of the present paper is to briefly survey recent work involving the authors of the present article and their collaborators that relates Pitman's transformation to classical integrable systems, such as the box-ball system, the ultra-discrete and discrete KdV equations, and the ultra-discrete and discrete Toda lattice equations \cite{CKST, CSMont, CS2, CS, CS3, CST, CST2}. Importantly, the framework that is developed across the latter articles enables the dynamics of the various systems to be initiated from infinite configurations, a novelty which is crucial for the study of invariant measures. In the special case of spatially independent and identically distributed (i.i.d.) configurations, the authors' progress upon the latter topic is reported here. Whilst the study of classical integrable systems from a probabilistic perspective has hitherto been somewhat limited, we note that beyond the results discussed, there have been increasing efforts in this direction, see \cite{FG, FG2, Ferrari, KL, KLO, Lev, LLPS} for some of the interesting developments.

Although in Pitman's original study the discrete-time version of Pitman's transformation was central to the proof of the main result concerning Brownian motion, the majority of the stochastic processes literature has focussed on the continuous-time setting. Part of the reason for this is that perhaps the most canonical discretisation of $T^{\int}$, whereby the integral is replaced by a sum, has somewhat different properties to the continuous-time version, cf.\ \cite{SS}. One of new aspects of our work is that through the connection to discrete integrable systems, we find an alternative, but still natural, discretisation of the operator $T^{\int}$ that has the appealing property of admitting an i.i.d.\ invariant measure for which the marginals of the process corresponding to $M-S$ precisely match those of the related continuous-time invariant measure (see Remark \ref{finrem}). In \cite{CS3}, we will discuss further advantages of this discretisation.

The remainder of the article is organised as follows. In Section \ref{pitsec} we present versions of Pitman's transformation that will appear later in the article. The corresponding discrete integrable systems are introduced in Section \ref{dissec}, and their relation with Pitman-type transformations is outlined in Section \ref{relsec}. Section \ref{invsec} incorporates a presentation of results on invariant measures for the systems under study, which are then summarised in Section \ref{sumsec}.

\section{Pitman-type transformations}\label{pitsec}

Similarly to the original version of Pitman's transformation that appears at \eqref{pitman}, the Pitman-type transformations that we consider will be of the form
\[T^*S=2M^*-S-2M^*_0\]
for some path-valued functional $M^*$. Note that, in contrast to the introduction, we will typically focus on two-sided paths $S:\mathbb{Z}\rightarrow\mathbb{R}$. Moreover, we will typically assume that $S_0=0$, and it is for this reason that we include the constant shift by $2M^*_0$, which ensures that we also have $(T^*S)_0=0$. Specifically, in the two-sided discrete-time case, namely for $S:\mathbb{Z}\rightarrow\mathbb{R}$, we will consider the following choices for $M^*$:
\begin{enumerate}
  \item[(a)]
  \[M(S)_n:=\sup_{m \le n}S_m;\]
  \item[(b)]
  \[M^{\vee}(S)_n:=\sup_{m \le n} \left(\frac{S_m+S_{m-1}}{2}\right);\]
  \item[(c)]
  \[M^{\sum}(S)_n:=\log \left(\sum_{m \le n} \exp\left(\frac{S_m+S_{m-1}}{2}\right)\right);\]
  \item[(d)]
  \[M^{\vee^*}(S)_{n}:=\left\{\begin{array}{ll}
  \sup_{m \le \frac{n-1}{2}}S_{2m}, & n\mbox{ odd}, \\
  \frac{M^{\vee^*}(S)_{n+1}+M^{\vee^*}(S)_{n-1}}{2}, &  n\mbox{ even};
  \end{array}
  \right.\]
  \item[(e)]
  \[M^{\sum^*}(S)_n:=\left\{\begin{array}{ll}
  \log \left(\sum_{m \le \frac{n-1}{2}} \exp\left(S_{2m}\right) \right), & n\mbox{ odd}, \\
  \frac{M^{\sum^*}(S)_{n+1}+M^{\sum^*}(S)_{n-1}}{2}, &  n\mbox{ even};
  \end{array}
  \right.\]
\end{enumerate}
writing the corresponding path transformations as $T^{\Z}$, $T^{\vee}$, $T^{\sum}$, $T^{\vee^*}$, $T^{\sum^*}$, respectively. Clearly, we need the right-hand sides in the definitions of $M^*$ to converge for these operators to be well-defined, and we note that this is always the case on the space of asymptotically linear functions, as given by
\begin{align}
\lefteqn{\mathcal{S}^{lin}:=}\label{slin}\\
&\left\{(S_n)_{n\in\mathbb{Z}}:\:\lim_{n\rightarrow-\infty}\frac{S_n}{n}\mbox{ and }\lim_{n\rightarrow+\infty}\frac{S_n}{n}\mbox{ exist and are strictly positive}\right\}.\nonumber
\end{align}
It will further be convenient to introduce shifted versions of $T^{\vee^*}$ and $T^{\sum^*}$, and so we define
\[\mathcal{T}^{\vee^*}(S):=\theta\circ{T}^{\vee^*}(S)-\theta\circ{T}^{\vee^*}(S)_0,\]
\[\mathcal{T}^{\sum^*}(S):=\theta\circ{T}^{\sum^*}(S)-\theta\circ{T}^{\sum^*}(S)_0,\]
where $\theta$ is the usual left-shift, i.e.\ if $x\in\mathbb{R}^\mathbb{Z}$, then $\theta(x)\in\mathbb{R}^\mathbb{Z}$ is given by $\theta(x)_n=x_{n+1}$. Again, we include a constant shift to ensure the new path passes through the origin.

\section{Discrete integrable systems}\label{dissec}

The discrete integrable systems that we consider in this article are based on discretisations (or ultra-discretisations) of the KdV and Toda lattice equations, as introduced in \cite{KdV} and \cite{Toda}, respectively. In each case, the variables in the model include those representing the configuration, $(z_n^t)_{n,t\in\mathbb{Z}}$ say, and those representing an auxiliary `carrier', $(w_n^t)_{n,t\in\mathbb{Z}}$ say, where the index $n$ refers to the spatial location of the variable, and the index $t$ refers to the temporal location. Moreover, the dynamics of the systems are locally-defined, in the sense that we have maps $(F_n)_{n\in\mathbb{Z}}$ and $m\in\mathbb{Z}$ such that $(z_{n-m}^{t+1},w_n^t)=F_n(z_n^t,w_{n-1}^t)$, or graphically, we have the following lattice structure:
\begin{equation}\label{lattice}
\hspace{-30pt}\boxed{F_n}\hspace{-15pt}\xymatrix@C-15pt@R-15pt{ & z_{n-m}^{t+1} & \\
            w_{n-1}^t \ar[rr]& & w_{n}^t.\\
             & z_n^t \ar[uu]&}
\end{equation}
We now describe a number of particular examples of such locally-defined dynamics, which, as we will explain in the next section, relate to the Pitman-type transformations of Section \ref{pitsec}. To keep the presentation concise, we refer the reader to \cite{CST2, IKT, TM, TT, TTMS, TH} for further background concerning the various models discussed.

\begin{enumerate}
\item[(a)] \emph{Box-ball system (BBS).} The most basic model in our discussion is the BBS, which was introduced in \cite{takahashi1990} as a simple example of a discrete model exhibiting the solitonic behaviour of the KdV equation. For this model, the configuration variables will be written $(\eta_n^t)_{n,t\in\mathbb{Z}}$. These take values in $\{0,1\}$, with $\eta_n^t=1$ representing the presence of a ball in the box at spatial location $n$ at time $t$, and $\eta_n^t=0$ representing the absence of such a ball. The carrier variables $(W_n^t)_{n,t\in\mathbb{Z}}$ are $\mathbb{Z}_+$-valued, with $W_n^t$ representing the number of balls transported from location $n$ to $n+1$ by the carrier between time $t$ and $t+1$. The evolution of the system is described as follows:
    \begin{equation}\label{bbsudkdv2}
    \begin{cases}
    \eta_n^{t+1}& = \min\{1-\eta_n^t,W^t_{n-1}\},\\
    W_n^t & =\eta_n^t+W^t_{n-1}-\eta_n^{t+1},
    \end{cases}
    \end{equation}
    which, roughly speaking, means that on each time step the carrier moves from left to right, picking up each ball it passes, and putting down a ball when it is carrying at least one and sees an empty box. There are many variations of the BBS, including the case of BBS($J$,$K$), where the boxes have capacity $J\in\mathbb{N}\cup\{\infty\}$ and the carrier has capacity $K\in\mathbb{N}\cup\{\infty\}$. Whilst we focus on the original model here (which can be thought of as BBS($1$,$\infty$)), it is in fact possible to describe the evolution of BBS($J$,$K$) in terms of a Pitman-type transformation whenever $J<K$ (see \cite{CS}). Other variations of the BBS incorporate balls of multiple colours. Again, such models can be described using Pitman-type transformations, but require higher dimensional versions to encode the relevant information \cite{Kondo}.
\item[(b)] \emph{Ultra-discrete KdV (udKdV) equation.} Generalising the BBS is the udKdV equation. In this model, both the configuration variables $(\eta_n^t)_{n,t\in\mathbb{Z}}$ and carrier variables $(U_n^t)_{n,t\in\mathbb{Z}}$ take values in $\mathbb{R}$, though one can also consider specialisations (including the BBS). Given a parameter $L\in\mathbb{R}$, the udKdV equation is given as follows:
    \begin{equation}\label{UDKDV}
    \begin{cases}
    \eta_n^{t+1}& = \min\{L-\eta_n^t,U_{n-1}^t\},\\
    U_n^{t} & =\eta_n^t+U_{n-1}^t-\eta_n^{t+1},
    \end{cases}
    \end{equation}
    which clearly matches \eqref{bbsudkdv2} when $L=1$. We can summarise the above equations by writing $(\eta^{t+1}_n,U_n^t)=F^{(L)}_{udK}(\eta^{t}_n,U^t_{n-1})$, which can be understood as locally-defined dynamics with a lattice structure for the variables as at \eqref{lattice} (with $m=0$). As with the BBS, there are other versions of the udKdV system that can be described by Pitman-type transformations; the $L$ parameter represents the box capacity, and we can also vary the carrier capacity, see \cite{CS3}.
\item[(c)] \emph{Discrete KdV (dKdV) equation.} Obtained from the KdV equation by a natural discretisation procedure, and yielding the udKdV equation through an ultra-discretisation procedure, is the dKdV equation. For this system, the configuration variables $(\omega_n^t)_{n,t\in\mathbb{Z}}$ and carrier variables $(U_n^t)_{n,t\in\mathbb{Z}}$ take values in $(0,\infty)$. Fixing $\delta\in(0,\infty)$, the dKdV equation is given by
    \[\begin{cases}
    \omega_n^{t+1} & =\left(\delta \omega_n^t +(U_{n-1}^t)^{-1}\right)^{-1},\\
    U_n^t & = U_{n-1}^t\omega_n^t (\omega_n^{t+1})^{-1}.
    \end{cases}\]
    Again, by writing $(\omega^{t+1}_n,U_n^t)=F^{(\delta)}_{dK}(\omega^{t}_n,U^t_{n-1})$, we can be understand the system being given by locally-defined dynamics with lattice structure as at \eqref{lattice} (with $m=0$). Moreover, variations of the dKdV equation with an additional parameter to above can also be handled within our framework, see \cite{CS3}.
\item[(d)] \emph{Ultra-discrete Toda lattice (udToda) equation.} Parallel to KdV-type equations, we consider Toda-type equations. The ultra-discrete Toda equation in particular has configuration variables $(Q_n^t,E_n^t)_{n,t\in\mathbb{Z}}$ and carrier variables $(U_n^t)_{n,t\in\mathbb{Z}}$ that take values in $\mathbb{R}$. The evolution of the system is described as follows:
    \[\begin{cases}
    Q_{n}^{t+1}=\min \{U_{n}^t,E_n^t\}, \\
    E_{n}^{t+1}=Q_{n+1}^t+E_{n}^t-Q_{n}^{t+1},\\
    U_{n+1}^t=U_{n}^t+Q_{n+1}^t-Q_{n}^{t+1},
    \end{cases}\]
    and summarised as $(Q_{n}^{t+1},E_n^{t+1},U_{n+1}^t)=F_{udT}(Q_{n+1}^t,E_n^t,U_{n}^t)$. We note that the udToda system also connects with the BBS in that we can view $Q_{n}^{t}$ as the length of the $n$th interval containing balls, and $E_n^t$ representing the length of the $n$th empty interval (at time $t$); of course, in the case of infinite balls, there is an issue of how to enumerate the intervals. For the udToda equation, the lattice structure is of the form
    \begin{equation}\label{qelattice}
    \xymatrix@C-15pt@R-15pt{ & Q_n^{t+1} & E_n^{t+1}& \\
            U_n^{t} \ar[rrr] & && U_{n+1}^t,\\
             & E_n^{t}\ar[uu]& Q_{n+1}^{t}\ar[uu]&}
    \end{equation}
    and so does not match \eqref{lattice}. However, it is possible to decompose the single map $F_{udT}$ with three inputs and three outputs into two maps $F_{udT^*}$ and $F_{udT^*}^{-1}$, each with two inputs and two outputs:
    \[\hspace{40pt}\xymatrix@C-15pt@R-15pt{\boxed{F_{udT^*}\vphantom{F_{udT^*}^{-1}}} & \min\{b,c\} &\boxed{F_{udT^*}^{-1}}& ^{a+b}_{-\min\{b,c\}}&\\
            c \ar[rr] && ^{c-\frac{b}{2}}_{-\frac{\min\{b,c\}}{2}}\ar[rr]&& ^{a+c}_{-\min\{b,c\}},\\
             & b \ar[uu]&&a\ar[uu]&}\]
    where we generically take $(a,b,c)=(Q_{n+1}^t,E_n^t,U_{n}^t)$. Including the additional lattice variables, we can thus view the system as locally-defined dynamics as at \eqref{lattice}, with $F_n$ alternating between $F_{udT^*}$ for $n$ even and $F_{udT^*}^{-1}$ for $n$ odd, and $m=1$.
\item[(e)] \emph{Discrete Toda lattice (dToda) equation.} Sitting between the original Toda lattice equation and its ultra-discrete version is the discrete Toda lattice equation, as given by:
    \[\begin{cases}
    I_n^{t+1}=J_n^t+U_n^t,\\
    J_n^{t+1}={I_{n+1}^{t}J_n^{t}}(I_n^{t+1})^{-1},\\
    U_{n+1}^t={I_{n+1}^{t}U_n^{t}}(I_n^{t+1})^{-1}.
    \end{cases}\]
    Here, the configuration variables $(I_n^t,J_n^t)_{n,t\in\mathbb{Z}}$ and carrier variables $(U_n^t)_{n,t\in\mathbb{Z}}$ take values in $(0,\infty)$, and we can summarise the above dynamics by $(I_{n}^{t+1},J_n^{t+1},U_{n+1}^t)=F_{dT}(I_{n+1}^t,J_n^t,U_{n}^t)$.
    Similarly to \eqref{qelattice}, in this case we have a lattice structure
    \[\xymatrix@C-15pt@R-15pt{ &I_n^{t+1}&J_n^{t+1} &\\
            U_n^{t} \ar[rrr]& & & U_{n+1}^{t},\\
     & J_n^t \ar[uu]&I_{n+1}^t\ar[uu]&}\]
     which can be decomposed into two maps, $F_{dT^*}$ and $F_{dT^*}^{-1}$, as follows:
     \[\xymatrix@C-15pt@R-15pt{\boxed{F_{dT^*}\vphantom{F_{udT^*}^{-1}}} & b+c &\boxed{F_{dT^*}^{-1}}& \frac{ab}{b+c}&\\
            c \ar[rr] &&\frac{c}{\sqrt{b^2+bc}} \ar[rr]&& \frac{ac}{b+c},\\
             & b \ar[uu]&&a\ar[uu]&}\]
     where we generically take $(a,b,c)=(I_{n+1}^t,J_n^t,U_{n}^t)$. So, again including the additional lattice variables, we can view the system as locally-defined dynamics as at \eqref{lattice}, with $F_n$ now alternating between $F_{dT^*}$ for $n$ even and $F_{dT^*}^{-1}$ for $n$ odd, and $m=1$.
\end{enumerate}

\section{Relation between locally-defined dynamics and Pitman-type transformations}\label{relsec}

For a system of locally-defined dynamics, as introduced at the start of the previous section, it is natural to consider an initial value problem of the following form: given initial condition $(z_n^0)_{n\in\mathbb{Z}}$, is it possible to find $(z_n^t)_{n,t\in\mathbb{Z}}$ and $(w_n^t)_{n,t\in\mathbb{Z}}$ such that $(z_{n-m}^{t+1},w_n^t)=F_n(z_n^t,w_{n-1}^t)$ holds for all $n,t\in\mathbb{Z}$? Moreover, if a solution exists, then is it unique?

Given the form of the equation $(z_{n-m}^{t+1},w_n^t)=F_n(z_n^t,w_{n-1}^t)$, it is clear that if both $(z_n^0)_{n\in\mathbb{Z}}$ and $(w_n^0)_{n\in\mathbb{Z}}$ are given, then we can compute $(z_n^1)_{n\in\mathbb{Z}}$ by setting $z_{n-m}^1=F^{(1)}_n(z_n^0,w_{n-1}^0)$, where the superscript $(1)$ refers to the first component of $F_n$. Hence, to solve the forward problem (where we only consider $t\geq 0$), it is enough to find a carrier $(w_n^0)_{n\in\mathbb{Z}}$ for $(z_n^0)_{n\in\mathbb{Z}}$ such that: $w_n^0=F^{(2)}_n(z_n^0,w_{n-1}^0)$ for all $n\in\mathbb{Z}$, and for which it is possible to repeat the procedure starting from the resulting $(z_n^1)_{n\in\mathbb{Z}}$. Moreover, in all the cases considered in this article, the locally-defined dynamics have a symmetry that means the backward problem can be solved in the same way.

Now, for any bijections $\mathcal{A}_n$ and $\mathcal{B}_n$ (defined on appropriate subsets of $\mathbb{R}$), we can rephrase the above problem in terms of the variables $x_n^t:=\mathcal{A}_n(z_n^t)$, $u_n^t:=\mathcal{B}_n(w_n^t)$. Indeed, for these variables, the equation of interest becomes
\[\left(x_{n-m}^{t+1},u_n^t\right)=K_n\left(x_n^t,u_{n-1}^t\right),\]
where $K_n := (\mathcal{A}_n \times \mathcal{B}_n) \circ F_n \circ (\mathcal{A}_n \times \mathcal{B}_{n-1})^{-1}$.
(The product map $\mathcal{A}_n \times \mathcal{B}_n$ is defined by setting $\mathcal{A}_n \times \mathcal{B}_n(a,b):=(\mathcal{A}_n(a),\mathcal{B}_n(b))$.) Of course, this is nothing but a change of variables. However, as a key assumption that enables a link to be made to a Pitman-type transformation, we will suppose that the variables have been changed in such a way that $K_n$ satisfies the conservation law:
\begin{equation}\label{cons}
K_n^{(1)}(a,b)-2K_n^{(2)}(a,b)=a-2b.
\end{equation}
In our work, we do not attempt to classify for which integrable systems such a change of variables exists, but we note that it is possible to do so for the systems considered here, i.e.\ BBS, udKdV, dKdV, udToda, dToda.

To relate to Pitman-type transformations, we need to introduce the path encoding of a configuration. In particular, given $(x_n)_{n\in\mathbb{Z}}$, we define the associated path $S:\mathbb{Z}\rightarrow \mathbb{R}$ by setting $S_0=0$ and
\[S_n-S_{n-1}=x_n,\qquad \forall n\in\mathbb{Z}.\]
Note that the mapping $x\mapsto S$ is one-to-one. By construction, it is easy to see that the existence of a carrier $(u_n)_{n\in\mathbb{Z}}$ for $(x_n)_{n\in\mathbb{Z}}$ (i.e.\ a solution to $u_n=K_n^{(2)}(x_n,u_{n-1})$) is equivalent to the existence of a path $M:\mathbb{Z}\rightarrow \mathbb{R}$ satisfying:
\begin{equation}\label{M-eq}
M_n=K^{(2)}_n\left(S_n-S_{n-1},M_{n-1}-S_{n-1}\right)+S_n.
\end{equation}
Indeed, if such an $M$ exists, then defining $(u_n)_{n\in\mathbb{Z}}$ by setting $u_n=M_n-S_n$, $n\in\mathbb{Z}$, gives a carrier for $(x_n)_{n\in\mathbb{Z}}$, and vice versa. Moreover, if such an $M$ exists and we define $(x^1_n)_n\in\mathbb{Z}$ to be the corresponding updated configuration (i.e.\ $x_{n-m}^1:=K_n^{(1)}(x_n,u_{n-1})$ for $u_n=M_n-S_n$), then the conservation law \eqref{cons} yields that
\begin{align*}
&x_{n-m}^{1}=2K_n^{(2)}\left(x_n,u_{n-1}\right)+x_n-2u_{n-1}\\
&=2K^{(2)}(S_n-S_{n-1},M_{n-1}-S_{n-1}) + (S_n-S_{n-1})-2(M_{n-1}-S_{n-1})\\
&=2M_n-S_n-(2M_{n-1}-S_{n-1}).
\end{align*}
Thus if we define $TS:\mathbb{Z}\rightarrow\mathbb{R}$ by setting $(TS)_n:=2M_n-S_n-2M_0$, then we see that the path encoding of $(x^1_n)_{n\in\mathbb{Z}}$ is precisely given by $\theta^m(TS)-\theta^m(TS)_0$, which is a (spatially-shifted) Pitman-type transformation of $S$ with respect to $M$.

To summarise, for any locally-defined dynamics $(K_n)_{n\in\mathbb{Z}}$ satisfying the conservation law at \eqref{cons}, we can associate a Pitman-type transformation of the corresponding path encoding. NB.\ By the change of variables, the conservation law can alternatively be expressed as
\[\mathcal{A}_n (F^{(1)}(a,b))-2\mathcal{B}_n (F^{(2)}(a,b))=\mathcal{A}_n(a)-2\mathcal{B}_n(b).\]
So far, this is simply a change of language. Importantly, however, in many examples, the equation \eqref{M-eq} can be solved explicitly. Moreover, it is often possible to determine uniquely a choice of $M$ for which the procedure can be iterated. In particular, in such cases, one obtains the existence and uniqueness of solutions to the initial value problem of interest.

\subsection{Examples}

It transpires that it is possible to follow the procedure described above in all of the examples of discrete integrable systems presented in Section \ref{dissec}. In particular, we are able to solve initial value problems for these systems whenever the path encoding of the initial configuration is an element of $\mathcal{S}^{lin}$, where the latter set was defined at \eqref{slin}. These results are presented in detail in \cite{CST2}. Rather than exhaustively repeat the exact statements here, however, we simply summarise how the framework applies in each case within the following table, and present one indicative example -- the udKdV equation -- in more detail (see Theorem \ref{udkdvthm} below). The maps $K^{\mathbb{Z}}$, $K^\vee$, $K^{\sum}$, $K^{\vee^*}$ and $K^{\sum^*}$ that appear in the table are given as follows: for $a,b\in\mathbb{R}$,
\begin{align*}
K^{\mathbb{Z}}(a,b)&:=\left(-\min\{a,2b-1\},b-\frac{a}{2}-\frac{\min\{a,2b-1\}}{2}\right),\\
K^{\vee}(a,b)&:=\left(-\min\{a,2b\},b-\frac{a}{2}-\frac{\min\{a,2b\}}{2}\right),\\
K^{\sum}(a,b)&:=\left(2\log(e^{-\frac{a}{2}}+e^{-b}),b-\frac{a}{2}+\log(e^{-\frac{a}{2}}+e^{-b})\right),\\
K^{\vee^*}(a,b)&:= \left(-\min\{a,b\},b-\frac{a}{2}-\frac{\min\{a,b\}}{2}\right),\\
K^{\sum^*}(a,b)&:= \left(\log(e^{-a}+e^{-b}),b-\frac{a}{2}+\frac{\log(e^{-a}+e^{-b})}{2}\right).
\end{align*}
Note that since we only apply $K^{\mathbb{Z}}$ to elements of $\{-1,1\} \times \{0,1,2,\dots\}$, it would be possible to consider other expressions that are equivalent on this subset of $\mathbb{R}^2$.
%DC - added the above sentence.
The latter four of the maps are obtained by a change of variables from $F_{udK}^{(L)}$, $F_{dK}^{(\delta)}$, $F_{udT^*}$ and $F_{dT^*}$, respectively. Observe that, unlike $F_{udK}^{(L)}$ and $F_{dK}^{(\delta)}$, the maps $K^\vee$ and $K^{\sum}$ are parameter free.

\begin{center}
\begin{adjustbox}{angle=90}
\begin{tabular}{c|c|c|c|c|c|c|c}
   \parbox{40pt}{\centering \emph{Integrable system}}&$z_n^t$&$w_n^t$&$\mathcal{A}_n(a)$&$\mathcal{B}_n(b)$&$m$&$K_n$& \parbox{66pt}{\centering \emph{Pitman-type transformation}}\\
    &&&&&&&\\
  \hline
   &&&&&&&\\
   BBS & $\eta_n^t$ & $W_n^t$ & $1-2a$ & $b$ & $0$ & $K^{\mathbb{Z}}$ & $T^{\mathbb{Z}}$  \\
  &&&&&&&\\
   \hline
   &&&&&&&\\
   udKdV & $\eta_n^t$ & $U_n^t$ & $L-2a$ & $b-\frac{L}{2}$ & $0$ & $K^{\vee}$ & $T^{\vee}$  \\
  &&&&&&&\\
  \hline
   &&&&&&&\\
   dKdV & $\omega_n^t$ & $U_n^t$ & $-\log\delta-2\log a$ & $\log b+\frac{\log \delta}{2}$ & $0$ & $K^{\sum}$ & $T^{\sum}$  \\
   &&&&&&&\\
   \hline
    &&&&&&&\\
    udToda & \parbox{54pt}{\centering $z_{2n-1}^t=Q_n^t$\\$z_{2n}^t=E_n^t$} &
   \parbox{98pt}{\centering $w_{2n-1}^t=U_n^t$\\$w_{2n}^t=F_{udT^*}^{(2)}(E_n^t,U_n^t)$}
   & $(-1)^na$ & $b$ & $1$ & \parbox{83pt}{\centering $K_{2n-1}=(K^{\vee^*})^{-1}$\\$K_{2n}=K^{\vee^*}$} & $\mathcal{T}^{\vee^*}$  \\
   &&&&&&&\\
  \hline
   &&&&&&&\\
  dToda & \parbox{54pt}{\centering $z_{2n-1}^t=I_n^t$\\$z_{2n}^t=J_n^t$} &
   \parbox{98pt}{\centering $w_{2n-1}^t=U_n^t$\\$w_{2n}^t=F_{dT^*}^{(2)}(J_n^t,U_n^t)$}
   & $(-1)^{n+1}\log a$ & $-\log b$ & $1$ & \parbox{83pt}{\centering $K_{2n-1}=(K^{\sum^*})^{-1}$\\$K_{2n}=K^{\sum^*}$} & $\mathcal{T}^{\sum^*}$
\end{tabular}
\end{adjustbox}
\end{center}
\renewcommand{\arraystretch}{1}

\newpage
As an illustrative example, in the next theorem, we present an application of the Pitman-type transformation approach to the udKdV equation.

\begin{thm} \label{udkdvthm}
Given $\eta=(\eta_n)_{n\in\mathbb{Z}}\in\mathbb{R}^\mathbb{Z}$, let $S$ be the path given by setting $S_0=0$ and $S_n-S_{n-1}=L-2\eta_n$ for $n\in\mathbb{Z}$. If $S\in\mathcal{S}^{lin}$, then there is a unique solution $(\eta_n^t,U_n^t)_{n,t\in\mathbb{Z}}$ to \eqref{UDKDV} that satisfies the initial condition $\eta^0=\eta$. This solution is given by
\[\eta_n^t:=\frac{L-S_n^t+S_{n-1}^{t}}{2},\qquad U_n^t:=M^\vee(S^t)_n-S_n^t+\frac{L}{2},\qquad \forall n,t\in\mathbb{Z},\]
where $S^t:=(T^\vee)^t(S)$ for all $t\in\mathbb{Z}$.
\end{thm}

\begin{rem}
Each of the operators $T^\vee$, $T^{\sum}$, $\mathcal{T}^{\vee^*}$ and $\mathcal{T}^{\sum^*}$ is a bijection from $\mathcal{S}^{lin}$ to itself. In particular, if $T^*$ is any of these operators and $S\in\mathcal{S}^{lin}$, then $(T^*)^t(S)$ is well-defined for all $t\in\mathbb{Z}$. Similarly, $T^\mathbb{Z}$ is a bijection on $\mathcal{S}^{lin}\cap \{S : \Z \to \Z\::\:S_0=0,\: |S_n-S_{n-1}|=1\}$.
\end{rem}

\begin{rem}
The statement of Theorem \ref{udkdvthm} is given in a slightly different form to the corresponding result in \cite{CST2}. This is because in the latter article, the definition of the transformation $T^*$ does not include the constant shift downwards by $2M^*_0$.
\end{rem}

\begin{rem}
As is discussed in detail in \cite{CST2}, the assumption that $S\in\mathcal{S}^{lin}$ includes various known solutions to the udKdV equation, such that as those based on finite or periodic initial configurations. The main novelty of our framework is that it also allows us to consider infinite initial configurations. As a particularly important example, one might consider $\eta=(\eta_n)_{n\in\mathbb{Z}}$ to be a random sequence that is stationary and ergodic under spatial shifts. If the density condition $\mathbf{E}\eta_0<\frac{L}{2}$ holds, then one almost-surely has that $S\in\mathcal{S}^{lin}$, and so the system can be started from the initial configuration $\eta$.
\end{rem}

\begin{rem}
In the case of the udKdV equation, writing the conservation law at \eqref{cons} in terms of the original lattice variables gives
\[\eta_n^t+U_{n-1}^t=\eta_n^{t+1}+U_{n}^t,\]
which can be interpreted transparently as conservation of mass by the system. The corresponding conservation law for the dKdV equation gives
\[\log\omega_n^t+\log U_{n-1}^t=\log \omega_n^{t+1}+\log U_{n}^t,\]
which can again be understood as conservation of mass. In the case of udToda, putting together the conservation laws for $K^{\vee*}$ and $(K^{\vee*})^{-1}$ in a way that eliminates the variable that was not part of the original lattice yields
\[E_n^t-Q_{n+1}^t-2U_n^t=E_n^{t+1}-Q_n^{t+1}-2U_{n+1}^t,\]
which is a combination of the conservation laws for the mass and length of the interval to which the local dynamics applies. An essentially similar argument yields the corresponding conservation rule for the dToda system. See \cite{CST2} for further discussion of this point.
\end{rem}

\begin{rem}
Using a different path encoding, it is shown in \cite{CST} that the dynamics of the ultra-discrete Toda equation can also be expressed in terms of $T^\mathbb{R}$.
\end{rem}

\section{Invariant measures}\label{invsec}

For the discrete integral systems of Section \ref{dissec}, it is a natural question to ask for which random initial configurations do the dynamics of the model leave the distribution of the configuration unchanged. Or, given that the mappings from the configurations of the discrete integrable systems to their path encodings are one-to-one, one might equivalently ask what distributions on paths are invariant under the associated Pitman-type transformations. In this section, we describe two approaches for verifying the invariance of measures, one of which is based on the path encoding viewpoint, and the other involves working directly with the configurations. Once we have done this, we describe some of our results concerning invariant measures for the discrete integrable systems of Section \ref{dissec} and Pitman-type transformations of Section \ref{pitsec}.

\subsection{Approaches for establishing invariance}

We call the first approach the `three conditions theorem', since it establishes a connection between the invariance of $S$ under a Pitman-type transformation and two natural symmetry properties. To describe these, the precise statement involves two reflection operators: a map $S\mapsto R(S)$, where
\[R(S)_n=-S_{-n}\]
(which corresponds to reflection of the configuration); and a map $W\mapsto\tilde{R}(W)$, where
\[\tilde{R}W_n:=W_{-n}\]
(which will be applied to the carrier $M(S)-S$). The following result initially appeared as \cite[Theorem 1.8]{CKST} in the context of the BBS, and has been generalised in \cite{CS3}. Since the latter version is stated in an abstract framework, one needs to make some further assumptions to capture the kind of operator to which the result applies. We will not detail these here, but simply note that these include a `reversibility' condition (REV), which ensures an appropriate interplay between the operators $R$ and $T$, and a `local time' condition (LOC), which enables the function $M(S)-M(S)_0$ to be recovered from $M(S)-S$. Both (REV) and (LOC) hold for all the discrete integrable systems of Section \ref{dissec}.

\begin{thm}[Three conditions theorem]\label{3ct}
Suppose that $T(S)=2M(S)-S-2M(S)_0$ is a Pitman-type transformation operator satisfying (REV) and (LOC). It is then the case that, for any probability measure supported on the domain of $T$, any two of the following conditions imply the third:\\
\[R(S) \buildrel{d}\over{=}S;\qquad  \tilde{R}(W)\buildrel{d}\over{=}W;\qquad T(S) \buildrel{d}\over{=}S,\]
where $W=M(S)-S$.
\end{thm}

\begin{rem}
If $S$ has i.i.d.\ increments, then $R(S) \buildrel{d}\over{=}S$ automatically holds, and so the invariance of $S$ under $T$ follows from the invariance of $W$ under $\tilde{R}$.
\end{rem}

For the second approach, which was developed for BBS($J$,$K$) in \cite{CS} and generalised in \cite{CS2}, we directly consider the evolution under a system of locally-defined dynamics, as introduced at the start of Section \ref{dissec}. We call the result a `detailed balance condition', since it is reminiscent of the corresponding result for determining invariant measures of Markov chains. Note that, for a measurable function $F$ and measure $\mu$ on the same space, we define $F(\mu):=\mu\circ F^{-1}$.

\begin{thm}[Detailed balance condition]\label{dbc}
Suppose $(F_n)_{n\in\mathbb{Z}}$ are all bijections. Moreover, let $\mathcal{C}\subseteq\mathbb{R}^\mathbb{Z}$ be such that for each $(z_n)_{n\in\mathbb{Z}}\in \mathcal{C}$, the system of locally-defined dynamics $(F_n)_{n\in\mathbb{Z}}$ has a unique solution $(z_n^t,w_n^t)_{n,t \in \Z}$ with $z^0=z$.\\
(a) Suppose $F_n \equiv F$, and $z=(z_n)_{n\in\mathbb{Z}}$ is an i.i.d.\ sequence supported on $\mathcal{C}$. It is then the case that $z$ is invariant under the dynamics if and only if
\[F(\mu \times \nu)=\mu \times \nu,\]
where $\mu$ is the distribution of $z_0^0$ and $\nu$ is the distribution of $w_0^0$.\\
(b) Suppose $F_{2n} \equiv F$ and $F_{2n+1} \equiv F^{-1}$, and $z=(z_n)_{n\in\mathbb{Z}}$ is an alternating i.i.d.\ sequence (i.e.\ the terms in the sequence are independent, with the odd terms and the even terms each being identically distributed) supported on $\mathcal{C}$.  It is then the case that $z$ is invariant under the dynamics if and only if
\[F(\mu \times \nu)=\tilde{\mu} \times \tilde{\nu},\]
where $\mu, \tilde{\mu}$ are the distributions of $z_0^0, z_1^0$, and $\nu,\tilde{\nu}$ are the distributions of $w_{-1}^0,w^0_0$, respectively.
\end{thm}

\begin{rem} Each of the two approaches has particular advantages and disadvantages, some of which we now briefly discuss.\\
(a) For the original Pitman transform $T^{\Z}$, if we restrict the state space to $\{S : \Z \to \Z\::\:S_0=0,\: |S_n-S_{n-1}|=1\}$ (i.e.\ the original BBS), then both Theorems \ref{3ct} and \ref{dbc} can be applied. On the other hand, if we consider more general increments for $S$, then the associated locally-defined dynamics are not given by bijections, and so we cannot apply the detailed balance condition.\\
(b) The three conditions theorem is applicable to both discrete-time random paths $S:\mathbb{Z}\rightarrow\mathbb{R}$ and continuous-time random paths $S:\mathbb{R}\rightarrow\mathbb{R}$. However, the detailed balance condition works only in the discrete-time case.\\
(c) The three conditions theorem only gives a sufficient condition for invariance, whereas the detailed balance condition gives an equivalent condition when it applies. Hence the latter result can be used to characterise all invariant measures of a Pitman-type transformations with i.i.d., or alternating i.i.d., increments.\\
(d) Unlike the detailed balance condition, the three conditions theorem does not depend on $z$ being an i.i.d.\ or alternating i.i.d.\ sequence. In fact, in \cite{CKST}, we apply the three conditions theorem to establish the invariance of some random configurations that are not i.i.d.
\end{rem}

\subsection{Invariant measures for discrete integrable systems}

Since the discrete integrable systems of Section \ref{dissec} have many conserved quantities, each of the models will admit a rich array of invariant measures (cf.\ \cite{CSMont, FG, Ferrari}). In the following four theorems, we give examples of invariant configurations based on i.i.d.\ or alternating i.i.d.\ sequences. Verification of the results via the detailed balance condition appears in \cite{CS2}. The latter article incorporates further discussion concerning the complete characterisation of such invariant measures. In the statements of the results, $Z$ represents a normalising constant, which will in general be different on each appearance.

\begin{thm} Suppose $\eta=(\eta_n)_{n\in\mathbb{Z}}$ is an i.i.d.\ sequence with marginals given by one of the following distributions:\\
(a) for some $\lambda>0$ and $c\in \R$ such that $c < \frac{L}{2}$,
\[\mathbf{P}\left(\eta_n \in dx\right)=\frac{1}{Z}e^{-\lambda x} \mathbf{1}_{[c,L-c]}(x)dx;\]
(b) for $h>0$ such that $L=\ell h$ for some $\ell \in \Z$, $\lambda >0$ and $ k \in \Z$ such that $k < \frac{\ell}{2}$,
\[\mathbf{P}\left(\eta_n = hm\right)=\frac{1}{Z}e^{-\lambda m} \mathbf{1}_{\{k, k+1,\dots, \ell -k\}}(m).\]
It is then the case that the path encoding of $\eta$ takes values in $\mathcal{S}^{lin}$ almost-surely, and the distribution of $\eta$ is invariant under the dynamics of the udKdV equation.
\end{thm}

\begin{thm}
Suppose $\omega=(\omega_n)_{n\in\mathbb{Z}}$ is an i.i.d.\ sequence with marginals given by the distribution: for some $\lambda >0$ and $c >0$,
\[\mathbf{P}\left(\omega_n \in dx\right)=\frac{1}{Z}e^{-cx^{-1}-c\delta x}x^{-\lambda-1}\mathbf{1}_{(0,\infty)}(x)dx.\]
It is then the case that the path encoding of $\omega$ takes values in $\mathcal{S}^{lin}$ almost-surely, and the distribution of $\omega$ is invariant under the dynamics of the dKdV equation.
\end{thm}

\begin{thm}
Suppose $(Q,E)=(Q_n,E_n)_{n\in\mathbb{Z}}$ is an i.i.d.\ sequence (with $Q_n$ independent of $E_n$ and) with marginals given by one of the following distributions:\\
(a) for some $0<\lambda_2<\lambda_1$ and $c \in \R$,
\[\mathbf{P}\left(Q_n \in dx\right)=\frac{1}{Z}e^{-\lambda_1 x} \mathbf{1}_{[c,\infty)}(x)dx,\]
\[\mathbf{P}\left(E_n \in dx\right)=\frac{1}{Z}e^{-\lambda_2 x} \mathbf{1}_{[c,\infty)}(x)dx;\]
(b) for some $0<\lambda_2<\lambda_1$, $h>0$ and $k \in \Z$,
\[\mathbf{P}\left(Q_n =mh \right)=\frac{1}{Z}e^{-\lambda_1 m} \mathbf{1}_{\{k,k+1,k+2,\dots\}}(m),\]
\[\mathbf{P}\left(E_n =mh \right)=\frac{1}{Z}\exp^{-\lambda_2 m} \mathbf{1}_{\{k,k+1,k+2,\dots\}}(m).\]
It is then the case that the path encoding of $(Q,E)$ takes values in $\mathcal{S}^{lin}$ almost-surely, and the distribution of $(Q,E)$ is invariant under the dynamics of the udToda equation.
\end{thm}

\begin{thm}
Suppose $(I,V)=(I_n,V_n)$ is an i.i.d.\ sequence (with $I_n$ independent of $V_n$ and) with marginals given by the following distribution: for some $0<\lambda_2<\lambda_1$ and $c >0$,
\[\mathbf{P}\left(I_n \in dx\right)=\frac{1}{Z}e^{-cx}x^{\lambda_1-1}\mathbf{1}_{(0,\infty)}(x)dx,\]
\[\mathbf{P}\left(V_n \in dx\right)=\frac{1}{Z}e^{-cx}x^{\lambda_2-1}\mathbf{1}_{(0,\infty)}(x)dx.\]
It is then the case that the path encoding of $(I,J)$ takes values in $\mathcal{S}^{lin}$ almost-surely, and the distribution of $(I,J)$ is invariant under the dynamics of the dToda equation.
\end{thm}

\subsection{Invariant measures for Pitman's type operators}

The following theorem transfers the results of the previous section to the setting of path encodings. In the statement, $S$ is either a simple random walk path (i.e.\ $S_0=0$, and $(S_n-S_{n-1})_{n\in\mathbb{Z}}$ is an i.i.d.\ sequence) or an alternating random walk path (i.e.\ $S_0=0$, and $(S_n-S_{n-1})_{n\in\mathbb{Z}}$ are independent, with the odd terms and the even terms each being identically distributed).

\begin{thm}
(a) If $S$ is a simple random walk such that
\[\mathbf{P}\left(S_n-S_{n-1}=x\right)=\left\{\begin{array}{ll}
                                                p, & \mbox{if }x=a,\\
                                                1-p-q, &  \mbox{if }x=0,\\
                                                q, & \mbox{if }x=-a,
                                              \end{array}\right.\]
for some $a >0$ and $0 \le q < p \le 1$ satisfying $p+q \le 1$, then $S$ is invariant under $T^{\Z}$.\\
(b) If $S$ is a simple random walk such that
\[\mathbf{P}\left(S_n-S_{n-1}  \in dx \right)=\frac{1}{Z} e^{\lambda x} \mathbf{1}_{[-a,a]}(x) dx\]
for some $a>0$ and $\lambda >0$, then $S$ is invariant under $T^{\vee}$. Moreover, if $S$ is a simple random walk such that either
\[\mathbf{P}\left(S_n-S_{n-1} =mh \right)=\frac{1}{Z} e^{\lambda m}\mathbf{1}_{\{-k,-k+1,\dots, k-1,k\}}(m),\]
or
\[\mathbf{P}\left(S_n-S_{n-1} =mh \right)=\frac{1}{Z} e^{\lambda m}\mathbf{1}_{\{-k+\frac{1}{2},-k+\frac{3}{2},\dots, k-\frac{1}{2}\}}(m),\]
for some $h >0$, $\lambda >0$ and $k \in \N$, then $S$ is invariant under $T^{\vee}$.\\
(c) If $S$ is a simple random walk such that
\[\mathbf{P}\left(S_n-S_{n-1}  \in dx \right)=\frac{1}{Z} e^{\lambda x-a \cosh\left(\frac{x}{2}\right)}\]
for some $a>0$ and $\lambda >0$, then $S$ is invariant under $T^{\sum}$.\\
(d) If $S$ is an alternating random walk such that
\[\mathbf{P}\left(S_{2n}-S_{2n-1} \in dx\right) =\frac{1}{Z}e^{-\lambda_1 x}\mathbf{1}_{(a,\infty)}(x)dx,\]
\[\mathbf{P}\left(S_{2n+1}-S_{2n} \in dx\right) =\frac{1}{Z}e^{\lambda_2 x}\mathbf{1}_{(-\infty,-a)}(x)dx,\]
for some $0< \lambda_1 < \lambda_2 $ and $a \in \R$, then $S$ is invariant under $\mathcal{T}^{\vee^*}$. Moreover, if $S$ is an alternating random walk such that
\[\mathbf{P}\left(S_{2n}-S_{2n-1} =mh \right)=\frac{1}{Z} e^{-\lambda_1 m}\mathbf{1}_{\{k,k+1,k+2\dots,\}}(m),\]
\[\mathbf{P}\left(S_{2n+1}-S_{2n} =mh \right)=\frac{1}{Z} e^{\lambda_2 m}\mathbf{1}_{\{\dots,-k-1,-k\}}(m),\]
for some $h >0$, $0< \lambda_1 < \lambda_2 $ and $k \in \Z$, then $S$ is invariant under $\mathcal{T}^{\vee^*}$.\\
(e) Let $S$ be an alternating random walk such that
\[\mathbf{P}\left(S_{2n}-S_{2n-1} \in dx\right) =\frac{1}{Z}e^{\lambda_1 x-ae^{x}}dx,\]
\[\mathbf{P}\left(S_{2n+1}-S_{2n} \in dx\right) =\frac{1}{Z}e^{-\lambda_2 x-ae^{-x}}dx,\]
for some $0<\lambda_1 <\lambda_2$ and $a>0$, then $S$ is invariant under $\mathcal{T}^{\sum^*}$.
\end{thm}

\begin{rem}
For (b)-(e), there are essentially two parameters that uniquely determine the drift and the variance of the increment. (Here, for (d) and (e), by increment we mean $S_{2n}-S_{2n-2}$, the distribution of which does not depend on the parameters $a$ or $k$.) For (a), since adding the delta measure at $0$ to the increment distribution does not affect to the invariance of $S$ under $T^{\Z}$, we have an extra parameter.
\end{rem}

\begin{rem}
For (a), (d) and (e), all invariant measures for the relevant Pitman-type transformations with i.i.d.\ (in the case of (a)) or alternating i.i.d.\ (in the case of (d) and (e)) increments are given by the previous theorem \cite{CS2}.
\end{rem}

\section{Summary}\label{sumsec}

Summarising the above, the following table outlines increment distributions of simple random walks/alternating random walks that are invariant under the Pitman-type transformations of Section \ref{pitsec}. In all cases, the process $W=M(S)-S$ is a reversible Markov chain, and it is possible to compute the distribution of $W_0$; this is also included in the table. For further details, see \cite{CS2,CS3}.

\begin{tabular}{c|c|c|c}
     \parbox{40pt}{\centering \emph{Operator}}
 & $S_n-S_{n-1}$ & $W_0$ &  \parbox{40pt}{\centering \emph{Integrable system}\\\hspace{10pt}} \\
  \hline
\parbox{40pt}{\hspace{10pt}\centering \\$T^{\Z}$\\\hspace{10pt}} &  \parbox{100pt}{\centering \hspace{10pt}\\Bernoulli (on $\{-1,1\}$), plus $\delta_0$ \\\hspace{10pt}}  & \parbox{60pt}{\centering \hspace{10pt}\\Geometric\\\hspace{10pt}} &\parbox{40pt}{\hspace{10pt}\centering \\BBS\\\hspace{10pt}}\\
\hline
\parbox{40pt}{\hspace{10pt}\centering \\$T^{\vee}$\\\hspace{10pt}} &  \parbox{100pt}{\centering \hspace{10pt}\\Truncated exponential or truncated geometric\\\hspace{10pt}}&\parbox{60pt}{\centering \hspace{10pt}\\Exponential or geometric\\\hspace{10pt}}& \parbox{40pt}{\hspace{10pt}\centering \\udKdV\\\hspace{10pt}}\\
 \hline
\parbox{40pt}{\hspace{10pt}\centering \\$T^{\sum}$\\\hspace{10pt}} &  \parbox{100pt}{\centering \hspace{10pt}\\Log of generalised inverse Gaussian\\\hspace{10pt}} & \parbox{60pt}{\centering \hspace{10pt}\\Log of inverse gamma\\\hspace{10pt}}&\parbox{40pt}{\hspace{10pt}\centering \\dKdV\\\hspace{10pt}}\\
\hline
\parbox{40pt}{\hspace{10pt}\centering \\$\mathcal{T}^{\vee^*}$\\\hspace{10pt}} & \parbox{100pt}{\centering\hspace{10pt}\\ Exponential or geometric, with alternating signs and parameters\\ (NB.\ $S_{2n}-S_{2n-2}$: asymmetric Laplace)\\\hspace{10pt}} & \parbox{60pt}{\centering \hspace{10pt}\\Exponential or geometric\\\hspace{10pt}}& \parbox{40pt}{\hspace{10pt}\centering \\udToda\\\hspace{10pt}}\\
\hline
\parbox{40pt}{\hspace{10pt}\centering \\$\mathcal{T}^{\sum^*}$} & \parbox{100pt}{\centering \hspace{10pt}\\Log of gamma, with alternating signs and parameters\\ (NB.\ $S_{2n}-S_{2n-2}$:\\log of beta)} & \parbox{60pt}{\centering \hspace{10pt}\\Log of inverse gamma} & \parbox{40pt}{\hspace{10pt}\centering \\dToda}
\end{tabular}
\bigskip

\begin{rem}\label{brem}
For any $c>0$ and $v \ge 0$, if $S:\mathbb{R}\rightarrow\mathbb{R}$ is defined by setting $S_x=vB_x+cx$, where $B=(B_x)_{x\in\mathbb{R}}$ is a two-sided standard Brownian motion, then the distribution of $S$ is invariant under both $T^{\R}$ and $T^{\int}$, see \cite{HW,ocyba}. In the former case, $M(S)-S$ is a reversible reflected Brownian motion with drift, which has exponential stationary distribution. In the latter, $M^{\int}(S)-S$ is a reversible diffusion, with stationary distribution being log of inverse gamma.
\end{rem}

\begin{rem}\label{finrem}
The parameters and distributions of the various models can be connected by ultra-discretisation, a correspondence between the KdV- and Toda-type systems, and continuous scaling limits \cite{CS2,CS3}. Indeed, taking an appropriate choice of parametrisation, it is possible to scale the invariant measures of $T^{\sum}$ to arrive at the Brownian invariant measure for $T^{\int}$ of Remark \ref{brem}. Moreover, this can be done in such a way that both the discrete and continuous models admit the same distribution for $W_0$.
\end{rem}

\section*{Acknowledgements}

The research of DC was supported by JSPS Grant-in-Aid for Scientific Research (C), 19K03540. The research of MS was supported by JSPS Grant-in-Aid for Scientific Research (B), 19H01792. We thank the referee for their careful reading of the paper, which resulted in a number of improvements being made.

\providecommand{\bysame}{\leavevmode\hbox to3em{\hrulefill}\thinspace}
\providecommand{\MR}{\relax\ifhmode\unskip\space\fi MR }
\providecommand{\MRhref}[2]{%
  \href{http://www.ams.org/mathscinet-getitem?mr=#1}{#2}
}
\providecommand{\href}[2]{#2}

\end{document}